\documentstyle{article}
\newtheorem{theorem}{Theorem}

\def\QED{\quad\blackslug\lower 8.5pt\null}

\begin{document}
\noindent
Proceedings of the 4th Congress of Geometry\\
Thessaloniki, 00--00 (1996)

\vspace*{0.2in}

\begin{center}
{\large \bf  ON A NORMALIZATION OF A GRASSMANN MANIFOLD} 
\vspace*{2mm}

 Maks A.  Akivis and Vladislav V. Goldberg
\end{center}


{\footnotesize 

\noindent 
{\bf Abstract.} On the Grassmann  manifold 
$G(m, n)$  of $m$-dimensional subspaces of an $n$-dimensional 
projective space $P^n$,  a certain supplementary 
construction called the  normalization is considered. 
By means of this 
normalization, one can construct  the structure of 
a Riemannian or semi-Riemannian manifold 
or an affine connection on $G(m, n)$.

\vspace*{2mm}



}

\setcounter{equation}{0}

{\bf 1.} 
Let $U$ be an open domain of the Grassmann  manifold $G(m, n)$   
of dimension $\rho = (m + 1)(n - m)$ coinciding with 
the dimension of  $G(m, n)$. This domain can coincide with 
the entire   manifold $G(m, n)$  or can be its proper subset. 
The domain $U$ is said to be {\em normalized} if to each 
of its  $m$-dimensional subspaces $p$ there corresponds a chosen 
subspace $p^*$ of dimension $n - m - 1$ in the projective space 
$P^n$, such that $p^*$ does not have common points with $p$. 
The  subspace $p^*$ is called the {\em normalizing subspace} for 
the subspace $p$. We will denote a  normalized domain $U 
\subseteq G (m, n)$ by $U^\nu$. 

Since the  subspace $p^*$ belongs to the Grassmannian 
$G (n - m - 1, n)$, a normalization of the manifold $G(m, n)$ 
is defined by a {\em normalizing mapping} 
\begin{equation}\label{eq:1}
\nu: G (m, n) \rightarrow G (n - m - 1, n)
\end{equation}
given in the domain $U \subseteq G (m, n)$ and having 
 a domain or a submanifold $U^*$ of the Grassmannian 
$G (n - m - 1, n)$ as its image. We assume that 
the  mapping $\nu$ is differentiable. 

 The number $r = \dim \; U^*$ 
coincides with the rank of the  mapping $\nu$. 
Since $\dim G (n - m - 1, n) = \rho = (m + 1)(n - m)$, 
we have $0 \leq r \leq \rho$. If $r = \rho$, then 
$U^*$ is an open domain of the manifold $G (n - m - 1, n)$. 
If $0 < r < \rho$, then $U^*$ is a proper submanifold of  
$G (n - m - 1, n)$. If $r = 0$, then $U^*$ consists 
of one fixed subspace $p^*$ of dimension $n - m - 1$ 
in the projective space $P^n$. 

If $r = \rho$, the normalization is called {\em nondegenerate}. 
In this case, there is a one-to-one differentiable 
correspondence between the domains $U$ and $U^*$. 
If $0 < r < \rho$, then the complete preimage $\nu^{-1} (p^*)$ 
of the normalizing subspace $p^*$ is a differentiable submanifold 
of dimension $\rho - r$ on the Grassmannian $G (m, n)$. 
If $r = 0$, then the complete preimage $\nu^{-1} (p^*)$ coincides 
with the entire domain $U$.

For $m = 0$, we arrive at the normalization 
of the  projective space $P^n$ considered in  [N 50],  \S 60. 

{\bf 2.} In this paper we will use the following index ranges:
$$
0 \leq \xi, \eta, \zeta \leq n; \;\;
0 \leq \alpha, \beta, \gamma, \delta, \epsilon \leq m;
\;\; m+1 \leq i, j, k, l \leq n.
$$

 Let us write the equations of the normalizing mapping 
$\nu$ using differential forms. To this end, with the pair of 
subspaces $p$ and $p^*$ we associate a family of point frames 
$\{A_\xi\}$ in such a way that $A_\alpha \in p$ and $A_i 
\in p^*$. For each frame of this family, we have 
\begin{equation}\label{eq:2}
d A_\alpha = \omega_\alpha^\beta A_\beta + \omega_\alpha^i A_i, 
\;\;\; d A_i = \omega_i^\alpha A_\alpha + \omega_i^j A_j. 
\end{equation}
The forms $\omega^\eta_\xi$ satisfy the structure equations of 
the projective space $P^n$: 
\begin{equation}\label{eq:3}
d \omega_\xi^\eta = \omega_\xi^\zeta \wedge \omega_\zeta^\eta.
\end{equation}

The 1-forms $\omega_\alpha^i $ 
are basis forms of the frame bundle associated with 
the Grassmannian $G (m, n)$. The 1-forms $\omega^\alpha_i$  
defining displacements of the subspace $p^*$ 
are expressed in terms of 
the basis forms $\omega^\alpha_i$ by the relations
\begin{equation}\label{eq:4}
 \omega_i^\alpha = \lambda_{ij}^{\alpha \beta} \omega_\beta^j. 
\end{equation}
These relations are differential equations of the normalizing 
mapping (1). The coefficients $\lambda_{ij}^{\alpha \beta}$ 
form a square matrix of order $\rho = (m+1)(n-m)$, 
whose rank $r$ is equal to the rank of the mapping $\nu$: 
$\mbox{rank} \; (\lambda_{ij}^{\alpha \beta}) = r$.

Next, taking exterior derivatives of equations (4) by means of 
(3) and applying Cartan's lemma to the exterior quadratic equations obtained, we find that 
\begin{equation}\label{eq:5}
 \nabla \lambda_{ij}^{\alpha \beta} 
= \lambda_{ijk}^{\alpha \beta \gamma}  \omega_\gamma^k, 
\end{equation}
where $\lambda_{ijk}^{\alpha \beta \gamma} = 
\lambda_{ikj}^{\alpha  \gamma \beta}$ and 
$\nabla \lambda_{ij}^{\alpha \beta} = d \lambda_{ij}^{\alpha \beta} - \lambda_{ik}^{\alpha \beta} 
\omega_j^k - \lambda_{kj}^{\alpha \beta} \omega_i^k 
+ \lambda_{ij}^{\alpha \gamma} \omega_\gamma^\beta 
+ \lambda_{ij}^{\gamma \beta} \omega_\gamma^\alpha$.
 If we fix an 
$m$--pair $(p, p^*)$, then equations (5) take the form 
$
 \nabla_\delta \lambda_{ij}^{\alpha \beta}  = 0, 
$
where $ \nabla_\delta \lambda_{ij}^{\alpha \beta} = 
 \nabla \lambda_{ij}^{\alpha \beta} (\delta)$ and 
$\delta = d_{\omega_\alpha^i = 0}$. The last relations 
show that the coefficients $\lambda_{ij}^{\alpha \beta}$ form a 
tensor connected with a first order 
differential neighborhood of the $m$-pair $(p, p^*)$. 
It is  called the {\em fundamental tensor of the 
normalized domain $U^\nu$}.

The object $\lambda_{ijk}^{\alpha \beta \gamma}$ occurring 
in equations (5) is also a tensor  connected with a 
second order differential neighborhood of the normalized 
Grassmann manifold.

{\bf 3.} In the  domain $U^\nu \subseteq G (m, n)$, 
we consider the quadratic differential form 
$
g = \omega_i^\alpha \omega_\alpha^i.
$
Substituting the values (4) of the forms $\omega^\alpha_i$ into 
the form $g$, we obtain
\begin{equation}\label{eq:6}
 g = g_{ij}^{\alpha \beta} 
\omega_\alpha^i \omega_\beta^j, 
\end{equation}
where the coefficients $g_{ij}^{\alpha \beta} $ 
are obtained if one symmetrizes the tensor 
$\lambda_{ij}^{\alpha \beta}$ 
simultaneously with respect to both lower and 
upper pairs of indices: 
$
g_{ij}^{\alpha \beta} 
= \frac{1}{2} (\lambda_{ij}^{\alpha \beta} 
+ \lambda_{ji}^{\beta \alpha}).  
$
Hence the quantities $g_{ij}^{\alpha \beta}$ 
themselves form a tensor which is symmetric with respect to 
 these pairs of indices. In view of this  the  
quadratic differential form $g$ is invariant in the 
domain $U^\nu$. 

Denote the rank of the matrix of coefficients of the quadratic 
form $g$ by $\widetilde{r}$. If $\widetilde{r} = \rho$, 
then the quadratic  form $g$ is nondegenerate and defines 
a Riemannian (or pseudo-Riemannian) metric in the domain $U^\nu$. 
If $\widetilde{r} < \rho$, then the  form $g$  defines 
a semi-Riemannian metric in the domain $U^\nu$ 
 for which the equation 
$
g_{ij}^{\alpha \beta}  \omega_\beta^j = 0
$
defines an isotropic distribution of dimension $\rho - 
\widetilde{r}$. 

The normalization $\nu$ is said to be {\em harmonic} 
if the coefficients in equations (4) are 
symmetric with respect to the vertical pairs of indices:
\begin{equation}\label{eq:7}
  \lambda_{ij}^{\alpha \beta} = \lambda_{ji}^{\beta \alpha}.
\end{equation}
If this is the case, then 
$g_{ij}^{\alpha \beta} 
= \lambda_{ij}^{\alpha \beta}$ and $\widetilde{r} = r$. 
If $r < \rho$, and the normalization $\nu$ is  harmonic, 
then the isotropic distribution defined by the  form $g$ 
is integrable, and its integral manifolds coincide with 
the complete preimages $\nu^{-1} (p^*)$ of the normalizing 
subspaces $p^*$.

{\bf 4.} Now we will establish a geometric meaning for 
the quadratic form (6). Consider the subspaces 
$p = A_0 \wedge A_1 \wedge \ldots \wedge A_m$ and 
$$
p' = (A_0 + dA_0) \wedge (A_1 + dA_1) \wedge \ldots 
\wedge (A_m + dA_m).
$$
 Their matrix coordinates $X$ and $Y$ (see  [R 96], Sect. 
{\bf 2.4.1}) are rectangular  matrices 
whose columns 
consist of the coordinates of points, determining 
the subspaces $p$ and $p'$, with respect to the 
frame ${\cal R} = \{A_0, A_1, \ldots , A_n\}$. 
These matrices are

\begin{equation}\label{eq:8}
X = \pmatrix{I_{m+1} \cr 
      O_{(n-m)\times (m+1)} \cr} \;\;\; \mbox{{\rm and}} \;\;\;
Y = \pmatrix{\delta_\beta^\alpha + \omega_\beta^\alpha  \cr 
    \omega_\beta^i \cr} \sim \pmatrix{\delta_\beta^\alpha    \cr 
         \omega_\beta^i \cr}, 
 \end{equation}
where $I_{m+1}$ is the identity matrix 
of order $m+1$, $O_{(n-m)\times (m+1)}$ is the rectangular zero  
$(n-m) \times (m+1)$ matrix, and the symbol $\sim$  denotes 
the equivalence of matrices with respect to multiplication from 
the right by a nondegenerate square matrix $(\delta_\beta^\alpha 
- \omega_\beta^\alpha)$ and discarding 
 second order terms with respect to the entries of the matrix 
$(\omega_\xi^\eta)$. 

Consider further the normalizing subspaces 
$p^* =  A_{m+1} \wedge \ldots \wedge A_n$
and
$$
p^{*\prime} =  (A_{m+1} + dA_{m+1}) \wedge \ldots \wedge 
(A_n + dA_n).
$$
It is easy to show that 
the tangential matrix coordinates $U$ and $V$ of $p$ and 
$p^{*'}$, 
 which are defined by the coefficients of linear equations 
which  the coordinates of points determining 
the subspaces $p*$ and $p^{*'}$ satisfy, can be reduced to the 
 forms:

\begin{equation}\label{eq:9}
U = \pmatrix{I_{m+1}  &  O_{(m+1) \times (n-m)} \cr} \;\;\;
\mbox{{\rm and}} \;\;\; 
V = \pmatrix{\delta_\beta^\alpha, - \omega_k^\alpha \cr}.
\end{equation}

In [R 96] (Sect. {\bf 2.4.4}) 
the cross-ratio $W$ of two $m$-pairs, whose  
matrix and tangential matrix coordinates are $X, Y$ and $U, V$, 
respectively,  was defined, and 
the following formula for its calculation was derived:
\begin{equation}\label{eq:10}
W = X (UX)^{-1} (UY) (VY)^{-1} V.
\end{equation}
 From the forms $X, Y, U$ and $V$, which we already calculated, 
we find that 
$$
U X = U Y = (\delta_\beta^\alpha), \;\;\; 
V Y = (\delta_\beta^\alpha - \omega_i^\alpha \omega_\beta^i), \;\; (V Y)^{-1} = (\delta_\beta^\alpha 
+ \omega_i^\alpha \omega_\beta^i).
$$
Thus, by applying formula (10), we find 
that  the cross-ratio $W$ of two  $m$-pairs 
$(p, p^*)$ and $(p', p^{*\prime})$ has the form:
\begin{equation}\label{eq:11}
W = \left(
\begin{array}{c}
\delta_\beta^\alpha + \omega_i^\alpha \omega_\beta^i \;\;\;
 - \omega_k^\alpha\\
O_{(n-m) \times (n+1)} 
\end{array}
\right).
\end{equation}
In the above calculations, we retain the terms of 
second order  with respect to the elements of the matrix 
$(\omega_\xi^\eta)$. Since such terms are principal, we 
discard the terms of order higher than two. 

To compute the quadratic form (6), we find the trace 
of the matrix $W$. It is: \\
$
\mbox{{\rm tr}}\;\; W = m + 1 +  \omega_i^\alpha \omega_\alpha^i. 
$
Since for small $x$ we have
$
\log (1 + x) \sim x,
$
it follows that
$$
\mbox{{\rm pr. p.}} \;
\log \Bigl(1 + \frac{1}{m+1} \omega_i^\alpha \omega_\alpha^i\Bigr) = 
\frac{1}{m+1}   \omega_i^\alpha 
\omega_\alpha^i,
$$
where $\mbox{{\rm pr. p.}}$ denotes the 
principal part of decomposition 
of the corresponding expression, and as a result, we find that
\begin{equation}\label{eq:12}
g = \omega_i^\alpha \omega_\alpha^i = (m + 1)\; 
\mbox{{\rm pr. p.}} \;
\Bigl[\log \Bigl(1 + \frac{1}{m+1} \mbox{{\rm tr}} \;\; W\Bigr)
\Bigr].
\end{equation}

Thus we have proved the following result.

\begin{theorem}
   The quadratic form $g$ is expressed 
in terms of the cross-ratio of two 
infinitesimally close $m$-pairs $(p, p^*)$ and 
$(p', p^{*\prime})$ by formula $(12)$. \rule{3mm}{3mm}
\end{theorem}

{\bf 5.} A normalization of the Grassmann  manifold $G(m, n)$ 
defines an affine connection on it. In fact, taking the exterior 
derivatives of  the basis forms  $\omega_\alpha^i$ of 
the manifold $G(m, n)$ and applying structure equations (3), 
we obtain 
\begin{equation}\label{eq:13}
d \omega^i_\alpha = \omega_\alpha^\beta \wedge 
\omega_\beta^i  + \omega_\alpha^j \wedge \omega_j^i 
= \omega_\beta^j \wedge (\delta_\alpha^\beta \omega_j^i   
- \delta_j^i  \omega_\alpha^\beta).
\end{equation}
Consider the 1-forms
\begin{equation}\label{eq:14}
 \omega^{i\beta}_{\alpha j} = \delta_\alpha^\beta \omega_j^i   
- \delta_j^i  \omega_\alpha^\beta.
\end{equation}
These forms are expressed in terms of the fiber forms 
$ \omega_\alpha^\beta$ and $\omega_j^i$   of the frame bundle 
associated with a domain $U^\nu \subseteq G (m, n)$. In the 
tangent space 
$T_p (\Omega)$, to the manifold $\Omega (m, n)$, which is the 
image of the   manifold $G(m, n)$ under the Grassmann mapping, 
these forms define a subgroup of the general linear 
group whose transformations preserve the 
cone of asymptotic directions of $G (m, n)$ determined by 
the equations $\omega_\alpha^i \omega_\beta^j - \omega_\alpha^j \omega_\beta^i = 0$.

Exterior differentiation of equations (14) leads to 
the  exterior  equations:
\begin{equation}\label{eq:15}
d \omega^{i\beta}_{\alpha j} - \delta_\alpha^\beta \omega_j^k   
\wedge \omega_k^i +  \delta_j^i  \omega_\alpha^\gamma 
\wedge \omega_\gamma^\beta  
 = \delta_\alpha^\beta \lambda_{jl}^{\gamma \epsilon}
 \omega_\epsilon^l \wedge \omega_\gamma^i -  \delta_j^i  \lambda_{kl}^{\beta \epsilon} \omega_\alpha^k 
\wedge \omega_\epsilon^l.
\end{equation}
The right-hand sides of equations 
(15) are expressed only in terms of 
the basis forms  $\omega_\alpha^i$ of the 
  domain $U^\nu$. By the facts 
from the  theory of spaces with affine connection 
(see, for example,  [KN 63], Ch. III), these 
equations show that the forms $\omega^{i\beta}_{\alpha j}$ 
define an affine connection on $U^\nu$, and the forms 
occurring in the right-hand sides of 
(15) are the {\em curvature forms} of this connection. 
Denote this   connection by $\Gamma^\nu$. 
The connection $\Gamma^\nu$ is uniquely determined by 
the normalization $\nu$. 
Note that affine connections on normalized Grassmannians were 
studied in [Ne 76].

Let us write the curvature forms of the connection $\Gamma^\nu$ 
in the form
\begin{equation}\label{eq:16}
\Omega^{i\beta}_{\alpha j} = (\delta_\alpha^\beta \delta_k^i 
 \lambda_{jl}^{\gamma \epsilon} + 
 \delta_\alpha^\gamma  \delta_j^i \lambda_{kl}^{\beta \epsilon}) 
 \omega_\epsilon^l \wedge \omega_\gamma^k.
\end{equation}
The alternated coefficients occurring in the right-hand sides of 
the last equations form the {\em curvature tensor} of the 
constructed connection. Equations (16) imply that this tensor 
has the following form: 
\begin{equation}\label{eq:17}
R_{\alpha jkl}^{i \beta \gamma \epsilon} 
= \frac{1}{2} \bigl(\delta_\alpha^\beta \delta_k^i 
 \lambda_{jl}^{\gamma \epsilon} + \delta_\alpha^\gamma \delta_j^i 
 \lambda_{kl}^{\beta \epsilon} 
- \delta_\alpha^\beta \delta_l^i 
 \lambda_{jk}^{\epsilon \gamma} - \delta_\alpha^\epsilon \delta_j^i  \lambda_{lk}^{\beta \gamma}\bigr),
\end{equation}
i.e. {\em this tensor is expressed only in terms 
of the components of the fundamental tensor of  the 
normalized  domain $U^\nu$}.

Equations (13) show that the affine connection $\Gamma^\nu$ 
is torsion-free. In view of this, the following theorem holds:

\begin{theorem}
The normalization $\nu$ of a normalized  domain $U^\nu 
\subseteq G (m, n)$  uniquely determines a torsion-free 
affine connection $\Gamma^\nu$ with the connection forms 
$(14)$ on it. The curvature tensor of this connection is linearly 
expressed in terms of the fundamental tensor of the normalization 
$\nu$ by formulas $(17)$. \rule{3mm}{3mm}
\end{theorem}

 Contracting the tensor (17) with respect 
to the indices $i, l$ and $\alpha, \epsilon$, we obtain 
the following expression for the Ricci tensor 
of the connection $\Gamma^\nu$:
\begin{equation}\label{eq:18}
R_{jk}^{\beta \gamma} = R_{\alpha jki}^{i \beta \gamma \alpha} 
= \frac{1}{2} \bigl(\lambda_{jk}^{\gamma \beta} 
+  \lambda_{kj}^{\beta \gamma} 
- (n + 1) \lambda_{jk}^{\beta \gamma}\bigr).
\end{equation}

Equations (18) imply the following result.

\begin{theorem}
The Ricci 
tensor of the connection $\Gamma^\nu$ is symmetric if and only if 
the normalization $\nu$ of the normalized domain $U^\nu$ 
is harmonic. \rule{3mm}{3mm}
\end{theorem}

{\bf 6.}  It is well-known that a Grassmann manifold 
$G (m, n)$ is a homogeneous space. However, in general, 
a normalized domain $U^\nu$ is not 
a homogeneous space. In fact, even two $m$-pairs 
$(p, p^*)$ and $(q, q^*)$ have a matrix invariant 
$W$---their cross-ratio. Thus, in general, there is no 
projective transformation superposing two neighborhoods 
$U (p, p^*)$ and $\widetilde{U}(\widetilde{p}, \widetilde{p}^*)$ 
of two $m$-pairs belonging to a normalized domain $U^\nu$. 

On the other hand, if a  normalized domain 
$U^\nu$ is homogeneous, then its fundamental 
tensor determining the location of an $m$-pair 
$(p', p^{* \prime})$, which is infinitesimally close to the 
$m$-pair $(p, p^*)$, must be covariantly constant, i.e. it must 
satisfy the condition
\begin{equation}\label{eq:19}
 \nabla \lambda_{ij}^{\alpha \beta} = 0,
\end{equation}
where $\nabla$ is the operator of  covariant differentiation 
with respect to the affine connection $\Gamma^\nu$.

Taking the 
exterior derivatives of the system of equations (19) 
by means of structure equations (3)  
and excluding the differentials 
$d  \lambda_{ij}^{\alpha \beta}$, we arrive at the 
system of relations:
\begin{equation}\label{eq:20}
\begin{array}{ll}
\lambda_{ik}^{\alpha \beta} \lambda_{jl}^{\gamma \epsilon} 
+ \lambda_{kj}^{\alpha \beta}  \lambda_{il}^{\gamma \epsilon}  
+ \lambda_{ij}^{\alpha \gamma} \lambda_{kl}^{\beta \epsilon} 
+ \lambda_{ij}^{\gamma \beta}  \lambda_{kl}^{\alpha \epsilon}\\
 - \lambda_{il}^{\alpha \beta} \lambda_{jk}^{\epsilon \gamma} 
- \lambda_{lj}^{\alpha \beta}  \lambda_{ik}^{\epsilon \gamma}    
- \lambda_{ij}^{\alpha \epsilon} \lambda_{lk}^{\beta \gamma} 
- \lambda_{ij}^{\epsilon \beta}  \lambda_{lk}^{\alpha \gamma} = 
0.    
\end{array}
\end{equation}

Thus, the following theorem is valid.

\begin{theorem}
   For the  normalization $\nu$ of the normalized 
domain $U^\nu$ with the fundamental tensor 
$\lambda_{ij}^{\alpha \beta}$ to be homogeneous it is  necessary 
and sufficient that the tensor 
$\lambda_{ij}^{\alpha \beta}$ satisfies the conditions $(19)$ and $(20)$. \rule{3mm}{3mm}
\end{theorem}

{\bf 7.} To find a solution of the system of 
equations (19) and (20), first we consider 
a   polar normalization, i.e. a normalization of 
the  Grassmann manifold $G (m, n)$ by means of a nondegenerate 
hyperquadric $Q$ of the space $P^n$ (see [N 50], \S\S 72--73).

Let $p_0$  be an $m$-dimensional subspace of  the space $P^n$ 
which is not tangent to  $Q$, and let 
$p_0^*$  be an $(n-m-1)$-dimensional subspace of  $P^n$ which is 
polar-conjugate to $p_0$ with respect $Q$. 
The subspaces $p_0$ and $p_0^*$ form a nondegenerate $m$-pair 
$(p_0, p_0^*)$. The set of  subspaces $p$, located in the 
same manner with respect to  $Q$  as $p_0$ 
(we will clarify  below the meaning of the expression ``in the 
same manner''), form an open domain $U$, and the subspaces $p^*$ 
polar-conjugate to the subspaces $p$ with  respect to  $Q$ 
define the {\em polar normalization} of this domain. 

If the  hyperquadric $Q$ is imaginary, then the domain $U$ 
coincides with the entire Grassmann manifold $G (m, n)$. 
Essentially, this case was studied in detail in  [L 61] 
where  the Riemannian geometry 
of the Grassmann manifold of subspaces of an Euclidean vector 
space was under investigation.

Let us associate a family of projective frames $\{A_\xi\}$ 
with an $m$-pair $(p, p^*)$ in such a way that 
the points $A_\alpha \in p$ and $A_i \in p^*$. We denote by 
$(A_\xi, A_\eta)$ the scalar 
product of the points $A_\xi$ and $A_\eta$ with 
respect to the  hyperquadric $Q$. 
Since the points $A_\alpha$ and $A_i$ are polar-conjugate with 
respect to this  hyperquadric, we have 
\begin{equation}\label{eq:21}
g_{i \alpha} = (A_i, A_\alpha) = 0.
\end{equation}
The scalar products
\begin{equation}\label{eq:22}
(A_i, A_j) = g_{ij} \;\; 
\mbox{{\rm and}} \;\; (A_\alpha, A_\beta) = g_{\alpha \beta}
\end{equation}
form nondegenerate symmetric matrices $(g_{\alpha \beta})$ and 
$(g_{ij})$. With respect to any  chosen frame, 
the equation of the  hyperquadric $Q$ can be written as 
\begin{equation}\label{eq:23}
g_{\alpha \beta} x^\alpha x^\beta + g_{ij} x^i x^j = 0.
\end{equation}
Moreover, the signature of each of the quadratic forms 
$g_{\alpha \beta} x^\alpha x^\beta$ and $g_{ij} x^i x^j$  
is not changed when the subspace $p$ moves in the 
normalized domain $U \subseteq G (m, n)$. 
This condition clarifies 
the meaning of the expression ``in the same manner'' which we 
used above to characterize the domain $U$. 


Differentiating equations (21) and (22) by 
means of equations (2), we find that 
$$
\omega_i^\alpha = - g^{\alpha \beta} g_{ij} \omega_\beta^j, \;\;
\nabla g_{ij} = 0, \;\;\;  \nabla g^{\alpha \beta} = 0,
$$
where $g^{\alpha\beta}$ is the inverse tensor 
of the tensor $g_{\alpha\beta}$. 
Comparing these with equations  (4), we obtain  
the fundamental tensor of the polar normalization:
\begin{equation}\label{eq:24}
\lambda_{ij}^{\alpha \beta} = - g^{\alpha \beta} g_{ij}. 
\end{equation}
Since the tensors $g^{\alpha \beta}$ and $g_{ij}$ are 
symmetric, this fundamental tensor satisfies condition 
(7), and {\em the polar normalization is harmonic}. 
Since the tensors $g^{\alpha \beta}$ and $g_{ij}$ are 
nondegenerate, {\em the fundamental tensor 
of the polar normalization is also nondegenerate}. 

 From relations (24) it follows that for the polar 
normalization we have
\begin{equation}\label{eq:25}
\nabla \lambda_{ij}^{\alpha \beta} = 0,
 \end{equation}
i.e. its  fundamental tensor is covariantly constant 
with respect to the connection $\Gamma^\nu$.  
Hence {\em the polar 
normalization of the  Grassmann manifold is homogeneous}. 

For the polar normalization, the  form (6) 
can be written as $g = - g^{\alpha \beta} g_{ij} \omega_\alpha^i 
\omega_\beta^j$. 
Thus, it is nondegenerate and defines a Riemannian (or 
pseudo-Riemannian) metric on the  normalized domain 
$U^\nu$ with a polar normalization $\nu$. By relation (25), 
{\em the connection $\Gamma^\nu$ is the Levi-Civita connection 
defined by this metric}.

Substituting values (24) of the fundamental tensor of 
the polar normalization into expressions (17), we obtain 
the following expression for the curvature tensor:
\begin{equation}\label{eq:26}
R_{ijkl}^{\alpha \beta \gamma \epsilon} =  \frac{1}{2} 
\bigl(g^{\alpha \beta} g^{\gamma \epsilon} (g_{il} g_{jk} 
- g_{ik} g_{jl}) 
 + (g^{\alpha  \epsilon} g^{\beta \gamma} 
- g^{\alpha \gamma} g^{\beta \epsilon}) g_{ij} g_{kl}\bigr).
\end{equation}
 Substituting  values (24) of the 
components of the fundamental tensor  of the polar 
normalization $\nu$ into (18), we find that 
$
R_{jk}^{\beta \gamma} = \frac{1}{2} (n - 1) 
g^{\beta \gamma}  g_{jk}, 
$
i.e. {\em the Ricci tensor of a polar-normalized domain 
 $U^\nu$ is proportional to its metric tensor}. 
But this means that {\em such a 
polar-normalized Grassmann manifold 
is an Einstein space}.

{\bf 8.} In conclusion, we consider the case 
when the normalizing mapping $\nu$ has zero rank: $r = 0$. 
Then the set of normalizing subspaces consists of 
a single subspace $p^*$ of dimension $n - m - 1$, and 
the normalized domain $U^\nu \subseteq G (m, n)$  consists of the 
$m$-dimensional subspaces $p$ not 
intersecting the normalizing subspace $p^*$. 

A projective space $P^n$, in which a subspace $p^*$ of dimension 
$n - m - 1$ is fixed, is called the {\em $m$-quasiaffine space}  
(see [R 59] and [D 88]) and is denoted by 
$A^n_m$.  The basis element of this space is 
a subspace $p$, and the entire space coincides with 
the domain $U^\nu$ considered above. 
The stationary subgroup of the element $p$ 
is  the  group  $H = {\bf GL} (m + 1) \times {\bf GL} (n - m)$. 
  
If we associate  a family of point frames with the subspace 
$p \in U^\nu$ in the manner indicated in Section {\bf 2}, 
then it is easy to prove that $\omega_i^\alpha = 0$, and thus 
$\lambda_{ij}^{\alpha\beta}=0$. It follows from (6) that $g = 0$, 
and {\em the  form $g$ 
 defines no metric in the domain $U^\nu$}. In addition,   
it follows from (17) that 
$R^{i \beta \gamma \epsilon}_{\alpha jkl} = 0$, and 
{\em   the connection $\Gamma^\nu$ is flat}. Thus, the domain 
$U^\nu$ is endowed with the structure of the affine space 
$A^\rho$ of dimension $\rho = (m+1)(n-m)$. But in this space 
the stationary subgroup $H$ leaves  invariant 
the Segre cone $SC_p$ with plane generators of dimensions $m + 1$ and $n - m$. This is the reason that this space is called the 
{\em Segre-affine space} and is denoted by $SA^\rho$.

Thus we have proved the following result.

\begin{theorem} Let $U^\nu$ be the domain of the Grassmann 
manifold $G (m, n)$ formed by its $m$-dimensional subspaces $p$ 
not having common points with a fixed subspace $p^*$ of dimension 
$n - m - 1$ $($the normalizing subspace$)$. Then the domain  
$U^\nu$ admits a mapping onto a Segre-affine space $S A^\rho$ 
which preserves the  structure of $U^\nu$. 
\rule{3mm}{3mm} 
\end{theorem} 

The mapping $s\colon U^\nu \rightarrow S A^\rho$ described 
in Theorem 5 is called the {\em stereographic projection} 
of the  Grassmann manifold $G (m, n)$.  
The stereographic projection of 
the Grassmann manifold $G (1, 3)$ was considered in 
 [SR 85], and for the general Grassmann 
manifold $G (m, n)$, it was considered in [S 32] 
(see also [D 88]). Since the Grassmann manifold 
$G (1, 3)$ is equivalent to the pseudoconformal space $C^4_2$, 
it admits the stereographic projection onto the pseudo-Euclidean 
space $R^4_2$ which is equivalent to the Segre-affine space 
$S A^4$.

\begin{center}
{\em Authors' addresses}:
\end{center}

{\sf             
M. A. Akivis, Department of Mathematics,            
Ben-Gurion University of the Negev,  
P.O. Box 653,  Beer-Sheva 84105, Israel       
}

{\em E-mail address}: akivis@black.bgu.ac.il

\vspace*{3mm}

{\sf 
V. V. Goldberg, 
 Department of Mathematics,   New Jersey Institute of Technology,  
 University Heights, Newark, NJ 07102, U. S. A.
}

{\em E-mail address}: vlgold@numerics.njit.edu

\end{document}